# Explicit solutions to the Navier-Stokes equation and the Euler equation

Yanyou Qiao

Aerospace Information Research Institute, Chinese Academy of Sciences,
Beijing 100101, China

qiaoyy@aircas.ac.cn

**Abstract:** In this paper, we obtain explicit solutions to the Navier-Stokes equation and the Euler equation. For any initial velocity $u_0(x)$ and the force vector $f(x,t)$, exact solutions can be explicitly solved as $u_i(x,t)=\sum_{n=0}^{\infty}u_{i,n}(x)t^n, (i=1,2,3)$, and $p(x,t)=\sum_{n=0}^{\infty}p_n(x)t^n$, where $u_{i,n}(x)\ (i=1,2,3;n=0,1,...)$ and $p_n(x),(n=0,1,...)$ are all known functions determined only by $u_0(x)$ and $f$. Explicit solutions of this paper can help in global regularity studies in the Navier-Stokes equation and the Euler equation. We admit that the convergence of these series still need to be proved before completely solving the existence of the Navier-Stokes equation and the Euler equation.

## 1. The Equations

The Euler and Navier–Stokes equations[1] describe the motion of a fluid in $R^3$. These equations are to be solved for an unknown velocity vector $u(x,t) = (u_1(x,t), u_2(x,t), u_3(x,t)) \in R^3$ and pressure $p(x,t) \in R$, defined for position $x = (x_1, x_2, x_3) \in R^3$ and time $t \geq 0$. For the incompressible fluids filling all of $R^3$, the Navier–Stokes equations are given by:

$$\frac{\partial u_i}{\partial t} + \sum_{j=1}^{3} u_j \frac{\partial u_i}{\partial x_j} = \upsilon \Delta u_i - \frac{\partial p}{\partial x_i} + f_i(x,t), (i = 1, 2, 3) \quad (1.1)$$

$$\frac{\partial u_1}{\partial x_1} + \frac{\partial u_2}{\partial x_2} + \frac{\partial u_3}{\partial x_3} = 0 \quad (1.2)$$

$$u(x,t)\,|_{t=0} = u_0(x) \quad (1.3)$$

Where $u_0(x) = (u_{1,0}(x), u_{2,0}(x), u_{3,0}(x))$ is a given, $C^\infty$ divergence–free vector field on $R^3$, it is the initial condition; $\Delta = \sum \frac{\partial^2}{\partial x_j^2}$ is the Laplacian in the space variables; $f = (f_1(x,t), f_2(x,t), f_3(x,t))$ is a given externally applied force; $\nu > 0$ is a positive coefficient (the viscosity), and the Euler equations are the above equations with $\nu = 0$.

For physically reasonable solutions, it is supposed that:

$$(p, u) \in C^\infty[R^3 \times [0, \infty)] \quad (1.4)$$

and $u(x,t)$ doesn't grow large as $|x| \to \infty$. Hence, the forces $f$ and the initial condition $u_0$ are restricted to satisfy

$$\left|\partial_x^\alpha u_0(x)\right| \leq C_{\alpha K}(1+|x|)^{-K} \quad (1.5)$$

on $R^3$ for any multi-index $\alpha$ and $K > 0$

$$\left|\partial_x^\alpha \partial_t^m f(x,t)\right| \leq C_{\alpha,m,K}\left(1+|x|+t\right)^{-K} \quad (1.6)$$

on $R^3 \times [0, \infty)$ for any multi-index $\alpha$, nonegative integer m, and $K > 0$

And there is a constant $C$ such that:

$$\int_{R^3} | u(x,t) |^2 dx < C \tag{1.7}$$

If there is a $(p, u)$ satisfying (1.1)-(1.7), we say that the Navier-Stokes equation has a regular solution; otherwise, there are blow-ups in the equation.

## 2. Preliminaries

The Schwartz space is define as

$$S = \{u(x) \in C^{\infty}(R^3) \mid \left|\partial_x^{\alpha} u(x)\right| \le C_{\alpha K}(1+ | x |)^{-K}, \forall \alpha, \forall K > 0\} \tag{2.1}$$

It can be seen that $u(x) \in S \Rightarrow \int_{R^3} | u | dx < \infty$ .

In this paper, we take $f_i(x,t) = \sum_{n=0}^{\infty} f_{i,n}(x) t^n, (i = 1, 2, 3)$ , where $f_{i,n} \in S$ .

It can be seen that $f$ is satisfied with (1.6).

The divergence free space is defined as

$$W = \{u(x) \in C^{\infty}(R^3) \mid \frac{\partial u_1}{\partial x_1} + \frac{\partial u_2}{\partial x_2} + \frac{\partial u_3}{\partial x_3} = 0\} \tag{2.2}$$

For any $u \in S$ , the projection of $u$ onto $W$ is defined as following

$$P(u) = u - \nabla(\omega) = \tilde{u} \tag{2.3}$$

Where $\omega \in C^{\infty}(R^3)$ is an unknown function to be determined. Since $div(\tilde{u}) = 0$ , the following poission equation can be got

$$\frac{\partial^2 \omega}{\partial^2 x_1} + \frac{\partial^2 \omega}{\partial^2 x_2} + \frac{\partial^2 \omega}{\partial^2 x_3} = \frac{\partial u_1}{\partial x_1} + \frac{\partial u_2}{\partial x_2} + \frac{\partial u_3}{\partial x_3} = -\varphi \tag{2.4}$$

And its solution is

$$\omega(x) = \int_{R^3} \frac{\varphi(\xi)}{4\pi | x - \xi |} d\xi \tag{2.5}$$

Thus $P(u) = \tilde{u} - \nabla(\omega)$ can be computed. It’s easy to check that

$$\int P(u).\nabla\omega dx = 0 \tag{2.6}$$

For any multi-index $\alpha$, there is

$$P(\partial_x^\alpha u) = \partial_x^\alpha (P(u)) \tag{2.7}$$

Lemma 2.1: ( By PG Lemarié-Rieusset)

suppose $\varphi(x) = O(1/|x|^4), (x \to \infty)$ and

$$\omega(x) = \int_{R^3} \frac{\varphi(\xi)}{4\pi |x-\xi|} d\xi$$

then $\partial_x \omega = O(1/|x|^2), (x \to \infty)$.

**Proof:** $\varphi(x) = O(1/|x|^4), (x \to \infty)$ is equal to

$\varphi(x) = O(1/(1+|x|)^4), (x \to \infty)$,hence, there is a constant C>0, such that

$\varphi \le C/(1/(1+|x|)^4)$.

We can see

$$\partial_{x_1}\omega(x) = -\frac{1}{4\pi}\int_{\mathbb{R}^3} \varphi(\xi)\frac{x_1-\xi_1}{|x-\xi|^3} d\xi$$

and $|\partial_{x_1}\omega(x)| \le \dfrac{C}{4\pi}\displaystyle\int_{\mathbb{R}^3} \frac{1}{(1+|\xi|)^4}\frac{1}{|x-\xi|^2} d\xi.$

Let $I = \displaystyle\int_{\mathbb{R}^3} \frac{1}{(1+|\xi|)^4}\frac{1}{|x-\xi|^2} d\xi$

Then $I = \displaystyle\int_{|x-\xi|\ge|\xi|} \frac{1}{(1+|\xi|)^4}\frac{1}{|x-\xi|^2} d\xi + \int_{|x-\xi|<|\xi|} \frac{1}{(1+|\xi|)^4}\frac{1}{|x-\xi|^2} d\xi = I_1 + I_2.$

If $|x-\xi| \ge |\xi|$

Then $|x| \le |\xi| + |x-\xi| \le 2|x-\xi|$

And thus $I_1 \le \dfrac{4}{|x|^2}\displaystyle\int_{\mathbb{R}^3} \frac{d\xi}{(1+|\xi|)^4}$

If $|x-\xi| < |\xi|$

Then $|x| \le |\xi| + |x-\xi| \le 2|\xi|$

And $\dfrac{1}{(1+|\xi|)^4}=\dfrac{1}{(1+|\xi|)}\dfrac{1}{(1+|\xi|)}\dfrac{1}{(1+|\xi|)^2}\leq 1\dfrac{2}{|x||\xi|}\dfrac{1}{|\xi|^2}$

And thus $I_2\leq\dfrac{2}{|x|}\int_{\mathbb{R}^3}\dfrac{1}{|x-\xi|^2}\dfrac{1}{|\xi|^2}d\xi$

with (writing $\xi=|x|R(\eta)$ and x=|x| R((1,0,0)) ,where R is a rotation)

$$\int_{\mathbb{R}^3}\frac{1}{|x-\xi|^2}\frac{1}{|\xi|^2}d\xi=\int_{\mathbb{R}^3}\frac{1}{|x|^2|(0,0,1)-\eta|^2}\frac{1}{|x|^2|\eta|^2}|x|^3\,d\eta=\frac{1}{|x|}\int_{\mathbb{R}^3}\frac{1}{|(0,0,1)-\eta|^2}\frac{1}{|\eta|^2}d\eta$$

Thus $I_2\leq\dfrac{C_1}{|x|^2}$

## 3. The series form solution

**Theorem 3.1**: There is a series form solution to the Navier-Stokes equation as following:

$$u_i(x,t)=\sum_{n=0}^{\infty}u_{i,n}(x)t^n,\ (i=1,2,3) \qquad (3.1)$$

$$p(x,t)=\sum_{n=0}^{\infty}p_n(x)t^n \qquad (3.2)$$

Where $\forall\alpha,\exists\partial_x^\alpha(u_{1,n},u_{2,n},u_{3,n})=1/|x|^2\ (|x|\to\infty),(n=0,1,2,...)$ , and

$(u_{1,n},u_{2,n},u_{3,n})\in C^\infty(R^3),(n=0,1,2,...)$.

**Proof:**

Let $u_i(x,t)=\sum_{n=0}^{\infty}u_{i,n}(x)t^n,(i=1,2,3)$, $p(x,t)=\sum_{n=0}^{\infty}p_n(x)t^n$ , we will solve the coefficients.

We will prove that all $u_{i,j}$ and $p_i$ can be got, and for any multi-index $\alpha$ , there is $\partial_x^\alpha u_{i,j}=O(1/|x|^2)\ (|x|\to\infty)$. This will be done by mathematical induction.

I) when i=0

Then $(u_{1,i}, u_{2,i}, u_{3,i}) = (u_{1,0}(x), u_{2,0}(x), u_{3,0}(x))$ is just the initial state. Since $(u_{1,0}, u_{2,0}, u_{3,0}) \in S$, it's clear that $\partial_x^\alpha (u_{1,0}, u_{2,0}, u_{3,0}) = O(1/\ |\ x\ |^2), (|\ x\ | \to \infty)$.

II) Suppose for all $0 \le i \le n-1$, $(u_{1,i}, u_{2,i}, u_{3,i})$ have been got, and $\partial_x^\alpha (u_{1,i}, u_{2,i}, u_{3,i}) = O(1/\ |\ x\ |^2)$, we want to solve $(u_{1,n}, u_{2,n}, u_{3,n})$ and $p_{n-1}$, and prove that $\partial_x^\alpha (u_{1,n}, u_{2,n}, u_{3,n}) = O(1/\ |\ x\ |^2)$.

By combining (1.1), (4.1) and (4.2), and compare the coefficients of $t^{n-1}$ in (1.1) and let:

$$v_{1,n} = \frac{1}{n}\begin{bmatrix} \upsilon \Delta u_{1,n-1} + f_{1,n-1} \\ -\sum_{i=0}^{n-1} u_{1,i} \dfrac{\partial u_{1,n-1-i}}{\partial x_1} - \sum_{i=0}^{n-1} u_{2,i} \dfrac{\partial u_{1,n-1-i}}{\partial x_2} - \sum_{i=0}^{n-1} u_{3,i} \dfrac{\partial u_{1,n-1-i}}{\partial x_3} \end{bmatrix} \tag{3.3}$$

$$v_{2,n} = \frac{1}{n}\begin{bmatrix} \upsilon \Delta u_{2,n-1} + f_{2,n-1} \\ -\sum_{i=0}^{n-1} u_{1,i} \dfrac{\partial u_{2,n-1-i}}{\partial x_1} - \sum_{i=0}^{n-1} u_{2,i} \dfrac{\partial u_{2,n-1-i}}{\partial x_2} - \sum_{i=0}^{n-1} u_{3,i} \dfrac{\partial u_{2,n-1-i}}{\partial x_3} \end{bmatrix} \tag{3.4}$$

$$v_{3,n} = \frac{1}{n}\begin{bmatrix} \upsilon \Delta u_{3,n-1} + f_{3,n-1} \\ -\sum_{i=0}^{n-1} u_{1,i} \dfrac{\partial u_{3,n-1-i}}{\partial x_1} - \sum_{i=0}^{n-1} u_{2,i} \dfrac{\partial u_{3,n-1-i}}{\partial x_2} - \sum_{i=0}^{n-1} u_{3,i} \dfrac{\partial u_{3,n-1-i}}{\partial x_3} \end{bmatrix} \tag{3.5}$$

Then

$$u_{1,n} = v_{1,n} - \frac{1}{n}\frac{\partial p_{n-1}}{\partial x_1}, \quad u_{2,n} = v_{2,n} - \frac{1}{n}\frac{\partial p_{n-1}}{\partial x_2}, \quad u_{3,n} = v_{3,n} - \frac{1}{n}\frac{\partial p_{n-1}}{\partial x_3} \tag{3.6}$$

By combining (1.2), (4.1) and (4.2), and compare the coefficients of $t^n$ in (1.2) ,we got:

$$\frac{\partial u_{1,n}}{\partial x_1} + \frac{\partial u_{2,n}}{\partial x_2} + \frac{\partial u_{3,n}}{\partial x_3} = 0 \tag{3.7}$$

By combining (3.6) and (3.7), we got the following poission equation:

$$\frac{\partial^2 p_{n-1}}{\partial x_1^2} + \frac{\partial^2 p_{n-1}}{\partial x_2^2} + \frac{\partial^2 p_{n-1}}{\partial x_3^2} = \varphi_{n-1} \tag{3.8}$$

Where

$$\begin{aligned}
&-\varphi_{n-1}(x) \\
&= \frac{\partial f_{1,n-1}}{\partial x_1} + \frac{\partial f_{2,n-1}}{\partial x_2} + \frac{\partial f_{3,n-1}}{\partial x_3} - \sum_{i=0}^{n-1} \frac{\partial u_{1,i}}{\partial x_1}\frac{\partial u_{1,n-1-i}}{\partial x_1} - \sum_{i=0}^{n-1} \frac{\partial u_{2,i}}{\partial x_2}\frac{\partial u_{2,n-1-i}}{\partial x_2} - \sum_{i=0}^{n-1} \frac{\partial u_{3,i}}{\partial x_3}\frac{\partial u_{3,n-1-i}}{\partial x_3} \\
&-2\sum_{i=0}^{n-1} \frac{\partial u_{2,i}}{\partial x_1}\frac{\partial u_{1,n-1-i}}{\partial x_2} - 2\sum_{i=0}^{n-1} \frac{\partial u_{3,i}}{\partial x_1}\frac{\partial u_{1,n-1-i}}{\partial x_3} - 2\sum_{i=0}^{n-1} \frac{\partial u_{2,i}}{\partial x_3}\frac{\partial u_{3,n-1-i}}{\partial x_2} \\
&= O(1/|x|^4), (|x| \to \infty)
\end{aligned}$$

Thus, the following integral can be got

$$p_{n-1}(x) = \int_{R^3} \frac{\varphi_{n-1}(\xi)}{4\pi |x-\xi|} d\xi \tag{3.9}$$

And

$$(u_{1,n}, u_{2,n}, u_{3,n}) = P((v_{1,n}, v_{2,n}, v_{3,n})) \tag{3.10}$$

From (3.6) and using lemma 2.1, there is $(u_{1,n}, u_{2,n}, u_{3,n}) = O(1/|x|^2)\ (|x| \to \infty)$.

Since $\partial_x^\alpha (u_{1,n}, u_{2,n}, u_{3,n}) = P(\partial_x^\alpha (v_{1,n}, v_{2,n}, v_{3,n}))$, there is also $\partial_x^\alpha (u_{1,n}, u_{2,n}, u_{3,n}) = O(1/|x|^2)\ (|x| \to \infty)$.

#

**Theorem 3.2**: Let $(\tilde{p}, \tilde{u})$ be the series defined in (3.1),(3.2):

$$\tilde{u}_i(x,t) = \sum_{n=0}^{\infty} u_{i,n}(x) t^n, (i = 1, 2, 3) \tag{3.11}$$

$$\tilde{p}(x,t) = \sum_{n=0}^{\infty} p_n(x) t^n \qquad (3.12)$$

$$R_0 = \{(x,t) \in R^3 \times [0,\infty) \mid |\tilde{u}| < \infty\} \qquad (3.13)$$

Then $(\tilde{p}, \tilde{u})$ is a Navier-Stokes fluid on $R_0$.

If $R_0 = \phi$, then $(\tilde{p}, \tilde{u})$ is the only regular solution.

**Proof:**

From the computation of all $u_{i,j}, (i = 1, 2, 3; j = 0, 1, \ldots)$ and

$p_i, (i = 0, 1, \ldots)$, it's clear that $(\tilde{p}, \tilde{u})$ is a Navier-Stokes fluid on $R_0$.

If $R_0 = \phi$, then $(\tilde{p}, \tilde{u})$ is a regular solution. By uniqueness analysis(see [2]), we have had the only regular solution $(p, u) = (\tilde{p}, \tilde{u})$, and it's in series forms.

#

However, the convergence of these series have not been proved.

Since the above process also hold when $\nu = 0$, a series form solution can also be got for the Euler equation, and their convergence also need further to be proved.

4. Discussions and conclusions

It can be seen that the most important point of this paper is the use of polynomial expansions along time t, so that the solution to the original complicated equations can be transformed into the solution of poisson equations whose solution is well-known. The polynomial expansions have been succeeded in getting a possible series solution to the Navier-Stokes equation and the Euler equation, but the convergence of these series have not been proved. The explicit solutions of this paper can help in global regularity studies in the Navier-Stokes equation and the Euler equation.

**Acknowledgement:** Thanks to my wife Honggan Wu and my son Rui Qiao for longtime full support; thanks to academician Guanhua Xu for continued inspiration in the research process; many thanks to professor PG Lemarié-Rieusset for helping a lot in this study.

References:

[1]. Charles L. Fefferman, EXISTENCE & SMOOTHNESS OF THE NAVIER–STOKES EQUATION, Official statement of the problem, Clay Mathematics Institute.

[2]. A. Bertozzi and A. Majda, Vorticity and Incompressible Flows, Cambridge U. Press, 2002.

[3]. Peter J. Olver, Introduction to Partial Differential Equations (Undergraduate Texts in Mathematics), springer, 2014.

[4]. L. Caffffarelli, R. Kohn, and L. Nirenberg, Partial regularity of suitable weak solutions of the Navier–Stokes equations, Comm. Pure & Appl. Math. 35 (1982), 771–831.

[5]. O. Ladyzhenskaya,, The Mathematical Theory of Viscous Incompressible Flows (2nd edition), Gordon and Breach, 1969.

[6]. F.–H. Lin, A new proof of the Caffffarelli–Kohn–Nirenberg theorem, Comm. Pure. & Appl. Math. 51 (1998), 241–257.

[7] A. Shnirelman, On the nonuniqueness of weak solutions of the Euler equation, Comm. Pure & Appl. Math.50 (1997), 1260–1286.